\documentclass[12pt]{amsart}
\usepackage{amssymb}
\usepackage{amsthm}
\usepackage{euscript}
\newcommand{\cA}{\EuScript{A}}
\newcommand{\cB}{\EuScript{B}}
\newcommand{\cC}{\EuScript{C}}

\newcommand{\cF}{\EuScript{F}}

\newcommand{\cL}{\EuScript{L}}
\newcommand{\cM}{\EuScript{M}}
\newcommand{\cN}{\EuScript{N}}

\newcommand{\cP}{\EuScript{P}}

\newcommand{\cR}{\EuScript{R}}
\newcommand{\cS}{\EuScript{S}}

\newcommand{\cV}{\EuScript{V}}
\newcommand{\cW}{\EuScript{W}}

\theoremstyle{plain}
\newtheorem{lem}[subsection]{Lemma}

\newtheorem{thm}[subsection]{Theorem}

\newtheorem{cor}[subsection]{Corollary}

\begin{document}
\title[Extension of a Theorem of Kaplansky]{An extension of a theorem of Kaplansky}
\author{Heydar Radjavi and Bamdad R. Yahaghi}

\address{Department of Pure Mathematics, University of Waterloo, Waterloo, Ontario, Canada N2L 3G1 \newline
\indent Department of Mathematics, Faculty of Sciences, University of Golestan, Gorgan 19395-5746, Iran}
\email{hradjavi@uwaterloo.ca, bamdad5@hotmail.com,\newline  bamdad@bamdadyahaghi.com}

\keywords{Semigroup, Division ring, Spectra, Irreducibility, Triangularizability.}
\subjclass{
15A30, 20M20}

\bibliographystyle{plain}

\begin{abstract}
A theorem of Kaplansky asserts that a semigroup of matrices with entries from a field whose members all have singleton spectra is triangularizable. Indeed, Kaplansky's Theorem unifies well-known theorems of Kolchin and  Levitzki on simultaneous triangularizability of semigroups of unipotent and nilpotent matrices, respectively.
First, we present a new and simple proof of Kaplansky's Theorem over fields of characteristic zero. Next, we show that this proof can be adjusted to show that the counterpart of Kolchin's Theorem over division rings of characteristic zero implies that of Kaplansky's Theorem over such division rings. Also, we give a generalization of Kaplansky's Theorem over general fields. We show that this extension of Kaplansky's Theorem holds over a division ring $\Delta$ provided the counterpart of Kaplansky's Theorem holds over $\Delta$.

\end{abstract}

\maketitle 

\bigskip

\begin{section}
{\bf Introduction}
\end{section}

\bigskip

A theorem of Kaplansky (see \cite[Theorem H on p. 137]{K} or \cite[Corollary 4.1.7]{RR1}) unifies two previous results: that of Levitzki, stating that a semigroup of nilpotent matrices is triangularizable (see \cite[Thoerem 35 on p. 135]{K} or \cite{L}, or \cite[Theorem 1.3]{Y2} for  a simple proof), and that of Kolchin deducing the same conclusion for a semigroup of unipotent matrices, i.e., those of the form $I+N$, where $I$ is the identity matrix and $N$ is nilpotent (see \cite{Ko} or \cite[Theorem C on p. 100]{K}). First, we present a new and simple proof of Kaplansky's Theorem over fields of characteristic zero. We show that this proof can be adjusted to show that the counterpart of Kolchin's Theorem over division rings of characteristic zero implies that of Kaplansky's Theorem over such division rings. 
Next, we give a generalization of Kaplansky's Theorem. To be more precise,  we prove that any semigroup of matrices with entries from a field of the form $T + N$, where $T$ comes from a triangularizable family $\mathcal{T}$ of matrices and $N$ is a nilpotent matrix coming from the commutant of $\mathcal{T}$ is triagularizable. This answers a question asked in \cite{RY} in the affirmative. Finally,  we show that our extension of Kaplansky's Theorem holds over a division ring $\Delta$ provided the counterpart of Kaplansky's Theorem holds over $\Delta$. 

Let us begin by fixing some standard notation. Let $\Delta$ be a division ring and $M_n(\Delta)$ the algebra of all $n \times n$ matrices over $\Delta$. The division ring  $\Delta$  could    in particular be a field. By a semigroup $\mathcal{S} \subseteq M_n(\Delta)$, we mean a set of matrices closed under multiplication. An ideal $\mathcal{J}$ of $\mathcal{S}$ is defined to be a subset of  $\mathcal{S}$ with the property that $SJ \in \mathcal{J}$  
and $JS \in \mathcal{J}$ for all $ S \in \mathcal{S}$ and $ J \in \mathcal{J}$. We view the members of $M_n(\Delta)$ as linear transformations acting on the left of $\Delta^n$, where $\Delta^n$ is the right vector space of all $n\times 1$ column vectors.  A semigroup $\mathcal{S}$ is called irreducible if the orbit of any nonzero $ x \in D^n$ under $ \mathcal{S}$ spans $\Delta^n$. When $n > 1$, this is equivalent to the members of $\mathcal{S}$, viewed as linear transformations on $\Delta^n$, having no common invariant subspace other than the trivial subspaces, namely, $\{0\}$ and $\Delta^n$. On the opposite of  irreducibility is triangularizability, when the common invariant subspaces of the members of $\mathcal{S}$ include a maximal subspace chain (of length $n$) in $\Delta^n$, i.e., there are subspaces 
$$\{0\} = \mathcal{V}_0 \subseteq  \mathcal{V}_1 \subseteq  \cdots \subseteq \mathcal{V}_n = \Delta^n,$$ 
where $ \mathcal{V}_j$ is a $j$-dimensional subspace invariant under every $S \in \mathcal{S}$. 
For a collection $\mathcal{C}$ in $M_n(D)$, by the commutant of $\mathcal{C}$, denoted by $\mathcal{C}'$, we mean 
$$\mathcal{C}' := \{A \in M_n(D) : AB = BA \   \forall \  B \in \mathcal{C}\}.$$

\bigskip 

We quote the following result  from \cite[Theorem 2.2.10]{Y3}  for reader's convenience. In fact the following theorem is a finite-dimensional version of \cite[Theorem 5]{Y} over general fields.

\bigskip 

\begin{thm} \label{1.1} 
{\it Let $\cV$ with $\dim \cV > 1$  be a finite-dimensional vector space over a field $F$  and $\cF$ a nonscalar triangularizable family of linear transformations on $\cV$. Then $\cF$ has a nontrivial hyperinvariant subspace.}
\end{thm}

\bigskip

\noindent {\bf Proof.} We note that for every family \(\cF\) of linear transformations 
\[{\cF}^{'}=({\rm Alg}({\cF}))^{'}.\] 
Thus \(\cF\) has a nontrivial hyperinvariant subspace iff \({\rm Alg}(\cF)\) does. Thus it suffices to prove the assertion for any  nonscalar triangularizable algebra, say \(\cA\), of linear transformations.
 
Now either the algebra  \(\cA\) is commutative or not. If it is a commutative algebra, note that by hypothesis there exists \( A \in \cA\) that is not scalar. Let \(\lambda\) be any eigenvalue of \(A\), and \(\cM\) the corresponding eigenspace of \(A\). Since \(\cA\) is commutative, for all \(B \in {\cA \cup \cA}^{'}\) and \(x \in \cM\) we have 
\[ABx=BAx={\lambda}Bx,\] 
i.e., \(Bx \in \cM\). So \(\cM\) is invariant under \({\cA \cup \cA}^{'}\). Now if the algebra \(\cA\) is not commutative, then there exist \(A, B \in \cA\) such that \(AB-BA \neq 0\). Set \(K_0=AB-BA\). Then \(K_0\) is a  nonzero nilpotent transformation in \(\cA\), for \(\cA\) is  triangularizable. 
Define  \({\cA}_1:={\cA}^{'}+ {\cA}*{\cA}^{'}\), where 
\[{\cA}*{\cA}^{'}:= \{ \sum_{i=1}^{k} A_iA_i' :k \in {\Bbb N} ,\ \ A_i \in {\cA},\ \  A_i' \in {\cA}^{'}, \ (1 \leq i \leq k) \}.\]
Clearly,  in view of the fact that \( {\cA}^{'}\) is a unital subalgebra of $\cL(\cV)$, we see that  \(\cA_1\) is a subalgebra of \( \cL(\cV)\) containing both \( \cA\) and \( {\cA}^{'}\). It thus suffices to prove that \(\cA_1\) has a nontrivial invariant subspace. 

We claim that ${\cA}_1K_0$, and hence ${\cA}_1K_0{\cA}_1$, the semigroup ideal generated by $K_0$ in \( {\cA}_1\), consists of nilpotents. To this end,  let  \( A_0= A' + \sum_{i=1}^{k} A_iA_i' \in {\cA}_1\) with 
\( A_i \in {\cA},\) where  \(  A', A_i' \in {\cA}^{'}, \ (1 \leq i \leq k, k \in {\Bbb N})\) be arbitrary. We prove that \(A_0K_0\) is nilpotent:
first of all we notice that  \( A_0K_0= A'K_0 + \sum_{i=1}^{k} A_{i_0}A_i'\), where \( A_{i_0}=A_iK_0 \in {\cA}\). 
Let $n= \dim \cV$. Set 
\[{\cS}:= \{A\in {\cA}: A^n=0\}.\]
Since $\cA$ is triangularizable, it follows that \(\cS\) is a nonzero semigroup ideal of  \(\cA\) consisting of nilpotent transformations (note that \(0 \not= K_0 \in \cS\)).
The set  \({\cS}{\cA}^{'}\) is indeed a semigroup consisting of nilpotents because for all \( A \in {\cA}, A' \in {\cA}^{'}\) we have \( AA'=A'A\) and that $\cS$ is a semigroup of nilpotents. Thus Levitzki's Theorem  shows that \({\cS}{\cA}^{'}\) is triangularizable. Therefore \({\rm Alg}({\cS}{\cA}^{'})\), the algebra generated by \({\cS}{\cA}^{'}\), consists of nilpotents.
 We have  
$$ A_0K_0= K_0A' + \sum_{i=1}^{k} A_{i_0}A_i',$$
 where \( A_{i_0}=A_iK_0 \in {\cA}\). In fact 
\( A_{i_0}=A_iK_0 \in {\cS}\), for \(K_0 \in \cS\) and \(\cA\) is triangularizable. Now clearly
\(A'K_0=K_0A' \in {\cS}{\cA}^{'}\) and \( A_{i_0}A_i'\in {\cS}{\cA}^{'}\). Therefore \( A_0K_0 \in {\rm Alg}({\cS}{\cA}^{'})\), and hence \(A_0K_0\) is a nilpotent transformation. Thus ${\cA}_1K_0{\cA}_1$ is a nonzero semigroup ideal of ${\cA}_1$ consisting of nilpotents which must be triangularizable, and hence reducible, by Levitzki's Theorem. Now reducibility of the nonzero ideal  ${\cA}_1K_0{\cA}_1$ implies that of ${\cA}_1$, completing the proof.
\hfill \qed 

\bigskip

The following is the counterpart of the preceding theorem over division rings.

\bigskip 

\begin{thm} \label{1.2} 
{\it Let $\cV$ with $\dim \cV > 1$  be a finite-dimensional vector space over a division ring $D$ with center $F$   and $\cF$ a triangularizable family of linear transformations on $\cV$ such that the $F$-algebra generated by $\cF$ contains a nonzero nilpotent linear transformation. Then $\cF$ has a nontrivial hyperinvariant subspace.}
\end{thm}

\bigskip

\noindent {\bf Proof.} The proof is an imitation of that of the preceding theorem, which is omitted for the sake of brevity. We refer the reader to \cite[Theorem 4.2.4]{Y3} for a detailed proof. \hfill \qed 

\bigskip 

Let $\cV$ be a finite-dimensional vector space over a division ring $D$.
For a triangularizable linear transformation $T \in  \cL(\cV)$, we say that $\lambda \in D$ is an {\it inner
eigenvalue of $T$ relative to a triangularizing
basis} $\cB$ for $T$ if $\lambda$ appears on the main diagonal of the matrix  of $T$ with respect to the basis $\cB$. It is easy to verify that if $\{S, T\} \subset  \cL(\cV)$ is triangularizable and $T$ and $S$ have inner-eigenvalues in the center of $D$, then $ST-TS$ is nilpotent. 

\bigskip

\begin{cor} \label{1.3} 
{\it Let $\cV$  with $\dim \cV > 1$ be a finite-dimensional vector space over a division ring $D$ with center $F$  and $\cF$ a nonscalar triangularizable family of linear transformations on $\cV$ with inner-eigenvalues in $F$. Then $\cF$ has a nontrivial hyperinvariant subspace.}
\end{cor}

\bigskip

\noindent {\bf Proof.} The assertion is easy if the family is commutative. If not, there exist $A, B \in \cF$ such that $C:= AB - BA \not= 0$. Then, clearly, $C $ is a nilpotent linear transformation and belongs to the $F$-algebra generated by $\cF$. The assertion now follows from Theorem \ref{1.2}. \hfill \qed 

\bigskip

A standard result in simultaneous triangularization over general fields is that the notion of triangularizability is preserved by passing to quotients. This result over division rings perhaps first appeared in \cite{S1}, but it is implicit in \cite[Lemma 1.5.2]{RR1} over fields and in \cite[Lemma 2.4]{S2} over division rings.  The proof given below is an imitation of the proof given  over fields in  \cite[Lemma 1.2.4]{Mom}.
Recall that if $\cV $ is a left (right)
vector space over $D$ and $\cM$ is a subspace of $\cV $, then 
$\cV /\cM:= \{ \ x + \cM : \  x \in \cV\}$  is called the {\it quotient space}.
If $A$ is a linear transformation on $\cV $ and $\cM \subset \cN$ are invariant
subspaces for $A$, then the quotient transformation $A_{\cN /\cM}$ on $\cN /\cM$
is defined by $A_{\cN /\cM} (x+ \cM) = Ax +\cM$ for each $x \in \cN $; the invariance of $\cM$ and $\cN$ under $A$ guarantees that $A_{\cN /\cM}$ is
well-defined. If $\cC$ is a collection of linear transformations
on $\cV $, and if $\cM \subset \cN$ are two invariant subspaces
for $\cC$, then the {\it collection of quotients of $\cC$ with
respect to} $\{ \cM,\cN\}$, denoted by $\cC_{\cN/\cM}$, is the set of all quotient
transformations $A_{\cN /\cM}$ on $\cN/\cM$, where $ A \in \cC$. We say that a property
$\cP$ is {\it inherited by quotients\/} if every collection of
quotients of a collection satisfying $\cP$ also satisfies $\cP$,
e.g., the properties nilpotency, commutativity, having rank $\leq
1$, etc are inherited by quotients. The following asserts that the property of triagularizability is inherited by quotients. 

\bigskip 

\begin{lem} \label{1.4} 
{\it Let $\cV$ be a finite-dimensional vector space over a division ring $D$, $\cC$ a triangularizable collection of linear transformations on $\cV$, and $\cM$ and $\cN$ with $\cM \subseteq \cN $ two invariant subspaces for the collection $\cC$. Then $\cC_{\cN/\cM}$ is triangularizable.}
\end{lem}

\bigskip

\noindent {\bf Proof.} Without loss of generality assume that $\dim \cV > 1$. Suppose $\cM$ and $\cN$ with $\cM \subseteq \cN $ and $ \dim \frac{\cN}{\cM} > 1$  are two invariant subspaces for the collection $\cC$. We need to show that there exists an invariant subspace $\cR$ of $\cC$ such that $ \cM \subsetneq \cR \subsetneq \cN$. To this end, let 
$$0 = \cV_0  \subsetneq \cV_1 \subsetneq \cdots \subsetneq \cV_{n-1} \subsetneq \cV_n = \cV,$$
where $ n = \dim \cV$, be a triangularizing chain for $\cC$. Set $\cW_i := \cN \cap (\cM + \cV_i)$, where $ 0 \leq i \leq n$. It is plain that each $\cW_i$ is an invariant subspace of $\cC$ and that 
$$\cM = \cW_0  \subseteq \cW_1 \subseteq \cdots \subseteq \cW_{n-1} \subseteq \cW_n = \cN.$$
Clearly, $ \dim \frac{\cW_i}{\cW_{i-1}} \leq  \dim \frac{\cV_i}{\cV_{i-1}} = 1$. This implies that there is an $ 1 < i < n$ such that $\cM  \subsetneq  \cR:=  \cW_i  \subsetneq \cN$.  This completes the proof. 
\hfill \qed 

\bigskip 

We need the following  useful lemma, which is a quick consequence of the preceding lemma, in the proof of one of our main results. 

\bigskip

\begin{lem} \label{1.5}
{\it Let $n \in \mathbb{N}$  and $\cF$ a family of block upper triangular matrices in $M_n(D)$. Then $\cF$ is triangularizable iff its diagonal blocks are triangularizable.}   
\end{lem}

\bigskip

\noindent {\bf Proof.} The proof, which is an quick consequence of Lemma \ref{1.4}, is omitted for the sake of brevity. 
\hfill \qed

\bigskip

\begin{section}
{\bf Main Results}
\end{section}

\bigskip

We start off with a simple proof of Kaplnasky's Theorem over  fields of zero characteristic. We recall that if a semigroup $\mathcal{S}$ of matrices is irreducible, then so is every nonzero ideal $ \mathcal{J}$ of $\mathcal{S}$ (see \cite[Lemma 2.1.10]{RR1}). 

\bigskip

\begin{thm} \label{2.1} {\bf (Kaplansky)}
Let $ n > 1$ and let $F$ be a field with characteristic zero and $ \mathcal S$ a semigroup in $ M_n(F)$ consisting of matrices with singleton spectra. Then the semigroup 
$ \mathcal S$ is triangularizable. 
\end{thm}

\bigskip

\noindent {\bf Proof.} By passing to $F^*\mathcal{S}$, where $F^* = F \setminus \{0\}$, we may assume that $\mathcal{S}$ is closed under scalar multiplications by the nonzero elements of $F$. We only need to show that $\mathcal S$ is reducible. If the semigroup $\mathcal{S}$ contains a nilpotent element, then reducibility of $\mathcal{S}$ follows from that of the nonzero semigroup ideal of  all nilpotent elements of  $\mathcal{S}$. So it remains to prove the assertion when $\mathcal{S}$ contains no nonzero nilpotent element.  It is then plain that $\mathcal{S}$ is reducible iff the set of all unipotent elements of $\mathcal{S}$ is reducible. Thus, in view of Kolchin's Theorem, we will be done as soon as we prove that the set of all unipotent elements of $\mathcal{S}$ forms a semigroup. To this end, let $ I + N_1 , I + N_2 \in \mathcal{S}$ be arbitrary unipotent elements. We can write $ (I+ N_1)(I + N_2) = cI + N' \in \mathcal{S}$, where $c \in   F^*$ and  $N' $ is a nilpotent matrix. We need to show that $c=1$. If $N_2 = 0$, we have nothing to prove. Suppose $N_2 \not= 0$ so that $N_2^k \not= 0$ but $N_2^{k+1} = 0$ for some $k \in \mathbb{N}$ with $k < n$. Thus
$$(I + N_2)^m = I + {m \choose 1} N_2 + \cdots + {m \choose k} N_2^k, $$
for all $ m \in \mathbb{N}$. Recall that ${m \choose k} := 0$ whenever $ m <k$. Clearly, we have $ (I+ N_1)(I + N_2)^m = c_mI + N_m' \in \mathcal{S}$  for all $ m \in \mathbb{N}$ with $c_m \in F$ an $n$th root of unity, i.e., $c_m^n = 1$ and $N'_m$ a nilpotent matrix. Since the set of the $n$th roots of unity in $F$ has at most  $n$ elements, we see that there exists a subsequence $(m_i)_{i=1}^\infty$  and an $l \in \mathbb{N}$ such that $ (I+ N_1)(I + N_2)^{m_i} = c_{m_l}I + N'_{m_l}$ for all $ i \in \mathbb{N}$.
Therefore, 
$$ \Big( (I+ N_1)(I + N_2)^{m} - c_{m_l}I\Big)^n = 0, $$
for infinitely many $ m \in \mathbb{N}$. Now, fix $ 1 \leq i , j \leq n$ and note that the $(i, j)$ entry of the matrix $ \Big( (I+ N_1)(I + N_2)^{m} - c_{m_l}I\Big)^n$ is a polynomial of degree $kn$ in $m$ having infinitely many roots, namely, $m_j$'s ($j \in \mathbb{N}$). This implies that   the $(i, j)$ entry of the matrix $ \Big( (I+ N_1)(I + N_2)^{m} - c_{m_l}I\Big)^n$ is zero for all $m \in \mathbb{N} \cup \{0\}$. Consequently, 
$$ \Big( (I+ N_1)(I + N_2)^{m} - c_{m_l}I\Big)^n = 0, $$
for all $ m \in \mathbb{N} \cup \{0\}$. Setting $m=0 ,1$ in the above, we obtain $ c = c_{m_l}=1$, which is what we want. 
\hfill \qed

\bigskip

\begin{thm} \label{2.2}
The counterpart of Kolchin's Theorem over division rings of characteristic zero implies that of Kaplansky's Theorem over such division rings. In other words, if every semigroup of unipotent matrices over a division ring  $\Delta$ of characteristic zero is triangularizable, then so is every semigroup of matrices of the form $ cI + N$, where $c$ comes from the center of $\Delta$ and $ N$ is a nilpotent matrix with entries from $\Delta$. 
\end{thm}

\bigskip

\noindent {\bf Proof.} Let $F$ denote the center of $\Delta$. Then, the proof is identical to that of the preceding theorem except that one should use the Gordon-Motzkin Theorem (\cite[Theorem 16.4]{La}) to get that 
$$ \Big( (I+ N_1)(I + N_2)^{m} - c_{m_l}I\Big)^n = 0, $$
for all $ m \in \mathbb{N} \cup \{0\}$. Also one must use the Dieudonn{\'e}  Determinant (see \cite[Corollary 20.1]{D}) to conclude that $c_m \in F$ is an $n$th root of unity whenever  $ (I+ N_1)(I + N_2)^m = c_mI + N' \in \mathcal{S}$  with $m \in \mathbb{N}$, $c_m \in F$, and $N_1, N_2$ and $N'$ nilpotent. 
\hfill \qed

\bigskip

Here is our extension of Kaplansky's Theorem. This theorem affirmatively answers a question raised in \cite{RY}.

\bigskip 

\begin{thm} \label{2.3}
Let $ n \in \mathbb{N}$ and let $F$ be a field and $\mathcal{T}$ a triangularizable set of  matrices in $ M_n(F)$,  $\mathcal{N}$ the set of all nilpotents in $ M_n(F)$, and $ \mathcal S$ a semigroup in $ M_n(F)$ consisting of matrices of the form $T + N$, where $ T \in \mathcal{T}$ and $ N \in \mathcal{T}' \cap \mathcal{N}$. Then the semigroup  $ \mathcal S$ is triangularizable. 
\end{thm}

\bigskip

\noindent {\bf Proof.} 
We view the elements of $M_n(F)$ as linear transformations on $F^n$ and proceed by induction on $n$, the dimension of the underlying space. The assertion trivially holds for $n=1$. Assume $ n > 1$ and that the assertion holds for such semigroups of linear transformations acting on spaces of dimension less than $n$. If $\mathcal{T}$ consists of scalar matrices, then $\mathcal{S}$ is triangularizable by Kaplansky's Theorem. If not, then by Theorem \ref{1.1}, $\mathcal{T}$ has a nontrivial hyperinvariant subspace.  Therefore, there exists a nontrivial direct sum decomposition  $F^n = \mathcal{V}_1 \oplus \mathcal{V}_2$ with respect to which 
$$T = \left(\begin{array}{cc} T_{11} & T_{12}\\ 0  & T_{22} \end{array} \right) , \  N = \left(\begin{array}{cc} N_{11} & N_{12}\\ 0  & N_{22} \end{array} \right) 
,$$
for all $ T \in \mathcal{T}$ and $ N \in \mathcal{T}' \cap \mathcal{N}$.
Thus, for each $S \in \mathcal{S}$,  with respect to the decomposition $F^n = \mathcal{V}_1 \oplus \mathcal{V}_2$, we can write
\begin{eqnarray*}
S &= & T_S + N_S\\
 & =&  \left(\begin{array}{cc} T_{11} & T_{12} \\ 0  & T_{22} \end{array} \right) + \left(\begin{array}{cc} N_{11} & N_{12}\\ 0  & N_{22} \end{array} \right) \\
 & =& \left(\begin{array}{cc} T_{11} + N_{11}& T_{12} + N_{12} \\ 0  & T_{22} + N_{22} \end{array} \right),
\end{eqnarray*}
 where 
$$T_S = \left(\begin{array}{cc} T_{11} & T_{12}\\ 0  & T_{22} \end{array} \right) , \  N_S = \left(\begin{array}{cc} N_{11} & N_{12}\\ 0  & N_{22} \end{array} \right).$$
 For $ 1 \leq i \leq 2$, use $\mathcal{S}_{ii}$, $\mathcal{T}_{ii}$, and $\mathcal{N}_{ii}$ to, respectively,  denote the set of all $(i, i) $ block entries of $S \in \mathcal{S}$, $T \in \mathcal{T}$, and $ N \in \mathcal{T}' \cap \mathcal{N}$ with $S = T+ N$. For $i = 1, 2$, let $ n_i = \dim \mathcal{V}_i$  so that $n = n_1 + n_2$ and $\mathcal{N}_{i} $ denote the set of all nilpotent linear transformations  on $\mathcal{V}_i$. Note that for each $i= 1, 2$, $\mathcal{T}_{ii}$ is triangularizable and $\mathcal{S}_{ii}$ is a semigroup of matrices of the form $ T_{ii} + N_{ii}$, where $T_{ii} \in \mathcal{T}_{ii}$, $N_{ii} \in  \mathcal{T}_{ii}' \cap \mathcal{N}_{i}$. By Lemma \ref{1.4}, $ \mathcal{S}$ is triangularizable iff  both $\mathcal{S}_{11}$ and $\mathcal{S}_{22}$ are triangularizable. But $\mathcal{S}_{ii}$ ($i=1, 2$) is triangularizable by the induction hypothesis because it consists of elements of the form $ T_{ii} + N_{ii}$, where $T_{ii} \in \mathcal{T}_{ii}$, $N_{ii} \in  \mathcal{T}_{ii}' \cap \mathcal{N}_{i}$,  $\mathcal{T}_{ii}$ is triangularizable, and $n_i < n$. This completes the proof. 
\hfill \qed

\bigskip

Here is what we can say over general division rings. 

\bigskip

\begin{thm} \label{2.4}
Let $ n \in \mathbb{N}$ and let $D$ be a division ring over which Kaplansky's Therem holds, $\mathcal{T}$ a triangularizable set of  matrices in $ M_n(D)$ with inner-eigenvalues in $F$, the center of $D$,  $\mathcal{N}$ the set of all nilpotents in $ M_n(D)$, and $ \mathcal S$ a  semigroup in $ M_n(D)$ consisting of matrices of the form $T + N$, where $ T \in \mathcal{T}$ and $ N \in \mathcal{T}' \cap \mathcal{N}$. Then the semigroup  $ \mathcal S$ is triangularizable. 
\end{thm}

\bigskip

\noindent {\bf Proof.} The proof, which we omit for the sake of brevity, is identical to that of the preceding theorem except that one has to make use of Corollary \ref{1.3} as opposed to Theorem \ref{1.1}.
\hfill \qed

\bigskip 

By a {\it Kaplansky semigroup} of matrices over a division ring $D$, we mean a semigroup of matrices of the form $cI + N$, where $c$ is in the center of $D$ and $N $  is a nilpotent matrix with entries from $D$. By \cite[Theorem 2.2]{RY}, every finite Kaplansky semigroup of matrices over a general division ring is triangularizable. This result together with the proof of the preceding theorem implies the following. 
\bigskip

\begin{thm} \label{2.5}
Let $ n \in \mathbb{N}$ and let $D$ be a division ring and $\mathcal{T}$ a triangularizable set of  matrices in $ M_n(D)$ with inner-eigenvalues in $F$, the center of $D$,  $\mathcal{N}$ the set of all nilpotents in $ M_n(D)$, and $ \mathcal S$ a finite semigroup in $ M_n(D)$ consisting of matrices of the form $T + N$, where $ T \in \mathcal{T}$ and $ N \in \mathcal{T}' \cap \mathcal{N}$. Then the semigroup  $ \mathcal S$ is triangularizable. 
\end{thm}

\bigskip

\noindent {\bf Proof.} In view of \cite[Theorem 2.2]{RY}, the assertion is a quick consequence of the proof of the preceding theorem. 
\hfill \qed

\end{document}